\documentstyle[amssymbols]{article}
\newtheorem{lemma}{Lemma}
\newcommand{\BbbC}{{\Bbb C}}
\newcommand{\BbbP}{{\Bbb P}}
\newcommand{\BbbZ}{{\Bbb Z}}

\newcommand{\mud}{\mu\left(d\right)}
\newcommand{\mueight}{\mu\left(8\right)}
\newcommand{\cosn}{\cos\left({2\pi\over n}\right)}
\newcommand{\sinn}{\sin\left({2\pi\over n}\right)}
\newcommand{\tC}{\tilde{C}}
\newcommand{\hesse}{\mathop{\rm hesse}\nolimits}
\begin{document}
\begin{center}
\title{A Projective Surface of Degree Eight\\
       with 168 Nodes}
\author{S.~Endra\ss\\
        Mathematisches Institut der Universit\"at Erlangen-N\"urnberg\\
        Bismarckstra\ss e 1 1/2, 91054 Erlangen, Germany}
\maketitle
\end{center}
\section*{Introduction}
Consider algebraic surfaces in complex projective
threespace $\BbbP_3$, denote by a {\em node} of such
a surface an ordinary double point and by $\mud$ the maximal number of
nodes of an algebraic surface of degree $d$ in $\BbbP_3$
with no further degeneracies.
This note shows that $\mueight\geq 168$ by giving an example of an octic
surface $X_8$ with $168$ nodes. $X_8$ is found within a seven parameter
family of $112$--nodal octic surfaces admitting dihedral symmetry of order
sixteen. This improves the estimate of $\mueight$ given by examples
of Gallarati (\cite{gallarati}, $\mueight\geq 160$) and Kreiss (\cite{kreiss},
$\mueight\geq 160$). On the other hand, using Miyaokas upper bound
 \cite{miyaoka} for the number of nodes of a projective surface, we get
$\mueight\leq 174$, thus
\[
    168\leq\mueight\leq 174.
\]
I have made excessive use of the computer algebra system
Maple V R3 computing $X_8$, in particular the calculations depend heavily
on some new resp.\ improved features of release 3
(\cite{maple}, pp.\ 22--24).
The construction of $X_8$ involves no free parameters, and in fact
D.\ van Straten calculated that $X_8$ is rigid using MACAULAY.

The octic $X_8$ is invariant under the group $D_8\times{\Bbb Z}_2$,
therefore invariant under the reflection group ${\Bbb Z}_2^3$ of
order eight. So $X_8$ is the eightfold cover of a projective
quartic surface with thirteen nodes. This construction is
referred to as the Segre--trick and it has been used extensively
to construct sextic surfaces with a given number of nodes, see
\cite{cataneseceresa} and \cite{barth}.

{}From the existence of $X_8$ one however cannot deduce the existence
of surfaces with $161$ up to $167$ nodes. It can be checked that the
family of surfaces $X_8$ is constructed from contains no such surfaces.
\medskip\\
{\em Acknowledgments: I would like to thank W.~Barth for
suggesting to me that surfaces with many nodes can be found among
$D_n$--symmetric surfaces.}
\section*{The Dihedral Groups}
Let the dihedral group $D_n$ of order $2n$ acting on $\BbbP_3$ be
generated by the rotation
\[
    \phi\colon\left(x:y:z:w\right)\longmapsto
    \left({\textstyle\cosn} x-{\textstyle\sinn} y:
          {\textstyle\sinn} x+{\textstyle\cosn} y:z:w\right)
\]
and the involution
\[
    \tau\colon\left(x:y:z:w\right)\longmapsto
    \left(x:-y:z:w\right).
\]
A surface $X=\left\{F=0\right\}\subset\BbbP_3$ will be called $D_n$--symmetric
if $F$ is $D_n$--invariant. The planes of reflection symmetry of
$D_n$ are exactly the $n$ planes
\[
    E_j=\left\{ \sin\left({{j\pi\over n}}\right)x=
                \cos\left({{j\pi\over n}}\right)y\right\},\qquad
    j=0,\ldots,n-1.
\]
If now $X$ is $D_n$--symmetric and for some
$j\in\left\{0,\ldots,n-1\right\}$ the plane curve $C=X\cap E_j$
has got a singularity in $p_0=\left(x_0:y_0:z_0:w_0\right)$,
then (assuming that after a rotation $E_j=\left\{y=0\right\}$,
thus $y_0=0$)
\[
    \frac{\partial F}{\partial\left\{x,z,w\right\}}
    \left(p_0\right)=
    \frac{\partial\left.F\right|_{E_j}}{\partial\left\{x,z,w\right\}}
    \left(p_0\right)=0.
\]
The reflection symmetry gives
$F\left(x,y,z,w\right)=F\left(x,-y,z,w\right)$, so
\[
    \frac{\partial F}{\partial y}\left(p_0\right)=
    -\frac{\partial F}{\partial y}\left(p_0\right)=0,
\]
therefore $p_0$ induces an orbit of singularities on $X$ with length
\[
    \left|\hbox{orbit}\left(D_n,p_0\right)\right|=\left\{
    \begin{array}{c@{\qquad}l}
        1 & \hbox{if\ }x_0=0, \\
        n/2 & \hbox{if $n$ is even and\ }z_0=w_0=0, \\
        n & \hbox{otherwise}.\\
    \end{array}\right.
\]
\section*{The Construction of $X_8$}
Let us begin with the seven--parameter family of octic surfaces, therefore
define $D_8$--invariant polynomials
\begin{eqnarray*}
    P & = & \prod_{j=0}^7\left(\cos\left(
            {{j\pi\over 4}}\right)x+
            \sin\left(
            {{j\pi\over 4}}\right)y-w\right)\\
      & = & {\frac{1}{4}\left ({x}^{2}-{w}^{2}\right )
            \left ({y}^{2}-{w}^{2}\right )\left (\left (x+y\right )^{2}-2
            \,{w}^{2}\right )\left (\left (x-y\right )^{2}-2\,
            {w}^{2}\right )}\\
    Q & = & \left (a\left (x^{2}+y^{2}\right )^{2}
            +\left (x^{2}+y^{2}\right )\left (b\,z^{2}+c\,zw+d\,w^{2}
            \right )\right.\\
      &   & \hspace{5ex}\left.+e\,z^{4}+f\,z^{3}w+g\,z^{2}w^{2}
            +h\,zw^{3}+i\,w^{4}\right )^{2}\\
\end{eqnarray*}
with parameters $a$, $b$, $c$, $d$, $e$, $f$, $g$, $h$, $i\in\BbbC$ and set
$F=P-Q$, $X=\left\{F=0\right\}$. $P$ vanishes exactly on all eight
planes $H_j=\left\{\cos\left({j\pi\over 4}\right)x
+\sin\left({j\pi\over 4}\right)y=w\right\}$, $j=0,\ldots,7$.
So $P$ vanishes to the second order
on the $28$ lines $H_j\cap H_k$, $0\leq j<k\leq 7$ and
$Q$ vanishes to the second order on a quartic surface.
Therefore for general values of $a,\ldots,i$ the polynomial $F$ vanishes
to the second order
on $4\cdot 28=112$ points. Hence $X$ has got $112$ nodes, all of them
lying on some symmetry plane $E_j$, $j\in\left\{0,\ldots,7\right\}$.
Because of symmetry
it is sufficient to consider only $E_0=\left\{y=0\right\}$ and
$E_1=\left\{x=\left(1+\sqrt{2}\right)y\right\}$. Substituting the equations
for $E_0$ and $E_1$ one gets homogeneous coordinates $\left(x:z:w\right)$
on $E_0$ and $\left(y:z:w\right)$ on $E_1$.

Now consider the two curves $C_j=X\cap E_j$, $j=0,1$. Then, as divisors,
we have:
\begin{eqnarray*}
    \left\{\left.P\right|_{E_0}=0\right\} & = &
L_1+L_2+2\left(L_3+L_4+L_5\right)\\
    L_{1/2} & = & \left\{ x=\pm w\right\} \\
    L_{3/4} & = & \left\{ x=\pm\sqrt{2}\,w\right\}\\
    L_5 & = & \left\{ w=0\right\}\\
\end{eqnarray*}
$C_0$ has got singularities in those twelve points where
$\left\{\left.Q\right|_{E_0}=0\right\}$ meets one of the lines $L_3$, $L_4$ or
$L_5$, so $C_0$ is a plane octic curve with twelve singularities
admitting reflection symmetry $\BbbZ_2$. Analogously,
\begin{eqnarray*}
    \left\{\left.P\right|_{E_1}=0\right\} & = & 2\left(M_1+M_2+M_3+M_4\right)\\
    M_{1/2} & = & \left\{y=\pm w\right\}\\
    M_{3/4} & = & \left\{y=\pm\left(\sqrt{2}-1\right)w\right\}\\
\end{eqnarray*}
so $C_1$ has got sixteen singularities and admits reflection symmetry
$\BbbZ_2$.

The first step is to set $c=f=h=0$, then $X$ is $D_8\times\BbbZ_2$--symmetric.
Now both $C_0$ and $C_1$ admit reflection symmetry $\BbbZ^2_2$ and
therefore can be constructed from plane quartic curves $\tilde{C}_0$ and
$\tilde{C}_1$ by applying the Segre--trick.  The second step is
to set $e=-1$ to mod out all projective transformations
$z\mapsto\lambda\cdot z$, $\lambda\in\BbbC^*$ from this family of octics.

Now call a singularity of $\tC_0$ or $\tC_1$ outside of the coordinate axes
a singularity in {\em general position}. Analogously call a point of contact
of $\tC_0$ or $\tC_1$ to one of the coordinate axes outside the
three points of intersection of those axes a point of contact in
{\em general position}. Applying
proposition 5 of \cite{cataneseceresa} gives:
\begin{itemize}
\item $\tC_0$ has got two nodes $s_1$ and $s_2$ and two contact points
      $t_1$ and $t_2$ to $\left\{w=0\right\}$ in general position.
\item $\tC_1$ has got four nodes $u_1,\ldots,u_4$ in general position,
      no three of them collinear and thus splits into two conics.
\item Every node in general position of $\tC_0$ or $\tC_1$
      induces an orbit of $16$ singularities on $X$.
\item Every point of contact in general position of either
      $\tC_0$ or $\tC_1$ to one of the
      coordinate axes $\left\{w=0\right\}$ or $\left\{z=0\right\}$
      induces an orbit of $8$ singularities on $X$.
\end{itemize}
All points $s_1$, $s_2$, $t_1$, $t_2$, and $u_1,\ldots,u_4$ induce
orbits of singularities on $X$, which will be denoted by the same letters.
All those points are singularities of $X$, so if the determinant
of the hesse matrix of $F$ in some point of such an orbit does not vanish,
all points in the corresponding orbit are nodes.
After substituting
\begin{eqnarray*}
        i & = & -{1\over 4}\left(8\left(3-2\,\sqrt{2}\right)
          \left(4\,a+b^2\right)
          +4\left(2-\sqrt{2}\right)
          \left(bg+2\,d\right)
            +g^2
            -i_1^2\right)\\
        d & = & -{1\over 32}\left( 256\,a+64\,b^2+16\,bg-\sqrt{2}\,d_1^2+
            \sqrt{2}\,i_1^2\right)\\
        a & = & -{1\over 64}\left( 16\,b^2+a_1^2-\left(2-\sqrt{2}\right)d_1^2
            -\left(2+\sqrt{2}\right)i_1^2\right)\\
\end{eqnarray*}
%
the determinants of the hesse matrices of the above points can be
computed (with Maple, of course):
\begin{eqnarray*}
    \det\left(\hesse\left(X,s_{1/2}\right)\right) & = &
        8\,a_1^2\left(4\,b+2\,g\pm a_1\right )\\
    \det\left(\hesse\left(X,t_{1/2}\right)\right) & = &
        const\cdot\left(
            a_1^2-\left(2-\sqrt{2}\right)d_1^2-\left(2+\sqrt{2}\right)i_1^2
        \right )\\
    &   &  \hspace{2ex}\left(
            4\,b\pm\sqrt{
            -a_1^2+\left(2-\sqrt{2}\right)d_1^2+\left(2+\sqrt{2}\right)i_1^2}
            \right)^{2}\\
    \det\left(\hesse\left(X,u_{1/2}\right)\right) & = &
        const\cdot\left(
            4\,b+\left(2+\sqrt{2}\right)g
            \pm\left(2+\sqrt{2}\right)i_1\right)i_1^2\\
    \det\left(\hesse\left(X,u_{3/4}\right)\right) & = &
        const\cdot\left(
            4\,b+\left(2-\sqrt{2}\right)g
             \pm\left(2-\sqrt{2}\right)d_1\right)d_1^2\\
\end{eqnarray*}
The third step is to set the remaining parameters to:
\[
    \begin{array}{lll}
        a=-{1\over  4}\left(1+  \sqrt{2}\right) &
        b= {1\over  2}\left(2+  \sqrt{2}\right) &
        d= {1\over  8}\left(2+ 7\sqrt{2}\right)\\[2ex]
        g= {1\over  2}\left(1- 2\sqrt{2}\right) &
        i=-{1\over 16}\left(1+12\sqrt{2}\right) \\
    \end{array}
\]
\smallskip\\
%
%
Now set $X_8=X$. It can be checked that
$\tC_0$ admits an additional node
$s_3=\left(8\left(\sqrt{2}-1\right):1:4\right)$
and an additional point of contact
$t_3=\left(1:0:2\right)$
to $\left\{z=0\right\}$,
whereas $\tC_1$ splits into one conic $K$ and two lines and therefore admits
one additional node
$u_5=\left( 2\left(3-2\,\sqrt{2}\right):3-2\,\sqrt{2}:4\right)$.
Moreover $K$ has points of contact to
$\left\{w=0\right\}$ in
$v_1=\left( 1:3+2\,\sqrt{2}:0\right)$
and to $\left\{z=0\right\}$ in
$v_2=\left( 1:0:4\right)$.
Computing determinants of hesse matrices results in:
\begin{eqnarray*}
    \det\left(\hesse\left(X,s_{1}\right)\right) & = & 128 \\
    \det\left(\hesse\left(X,s_{2}\right)\right) & = & 1152 \\
    \det\left(\hesse\left(X,s_{3}\right)\right) & = &
        -128\left(239-169\,\sqrt{2}\right)\\
    \det\left(\hesse\left(X,t_{1}\right)\right) & = & {1\over 4}\\
    \det\left(\hesse\left(X,t_{2}\right)\right) & = &
        {3\over 4}\left( 3+2\,\sqrt{2}\right)\\
    \det\left(\hesse\left(X,t_{3}\right)\right) & = & {\frac {9}{512}}\\
    \det\left(\hesse\left(X,u_{1}\right)\right) & = &
        512\left(1451+1026\,\sqrt{2}\right)\\
    \det\left(\hesse\left(X,u_{2}\right)\right) & = &
        512\left(99-70\,\sqrt{2}\right)\\
    \det\left(\hesse\left(X,u_{3}\right)\right) & = &
        512\left(11243+7950\,\sqrt{2}\right)\\
    \det\left(\hesse\left(X,u_{4}\right)\right) & = &
        4608\left(331+234\,\sqrt{2}\right)\\
    \det\left(\hesse\left(X,u_{5}\right)\right) & = &
        2\left(1451-1026\,\sqrt{2}\right)\\
    \det\left(\hesse\left(X,v_{1}\right)\right) & = &
        512\left(3363+2378\,\sqrt{2}\right)\\
    \det\left(\hesse\left(X,v_{2}\right)\right) & = &
        {1\over 128}\left( 2979+2106\,\sqrt {2}\right)\\
\end{eqnarray*}
This gives $16+16+8+8+8=56$ additional nodes, thus $168$ nodes
altogether.
\section*{$X_8$ is smooth away from the 168 nodes}
Now one checks (again with Maple) that the 168 nodes are the only ones
to appear on the eight planes $E_j$, $j=0,\ldots,7$:
\begin{itemize}
\item points of contact to coordinate axes can be computed explicitly,
\item the fact that $K$ is non degenerate can be
      checked by computing its determinant,
\item the fact that $\tC_0$ is irreducible can be checked by
      projecting through a node onto some line not containing this
      node and computing the branch points with multiplicities.
\end{itemize}
One also checks that the points of intersection of
$X_8$ with the line $L=\left\{\left(x:y:0:0\right)\mid
\left(x:y\right)\in\BbbP_1\right\}$ are smooth.
If $X_8$ would have a singularity $p_0$ outside of all eight planes $E_j$,
$j=0,\ldots,7$ and outside the line  $L$, then $p_0$ would induce
an orbit of sixteen singularities of $X_8$. Now consider the following lemma:
\begin{lemma}
    Let $Y=\left\{G=0\right\}$ be a $D_n$--symmetric surface of degree
    $n\geq 1$. Then $Y$ has got no orbit of $2n$ nodes.
\end{lemma}
\noindent{\bf Proof:}\
If $n\leq 3$ we have $\mu\left(n\right)<2n$, so let $n\geq 4$.
Assume $O_{2n}$ is an orbit of $2n$ nodes of $Y$. Then $O_{2n}$ is
contained in some plane curve $K$ of degree $\leq 2$.
Let $E$ be the plane containing $K$. After a projective transformation
we may assume $E=\left\{z=0\right\}$. Then $C=Y\cap E$ is a plane
curve of degree $n$ with $C.K\geq 4n>2n$, so $C=K+C'$. But
$C'.K\geq 2n>2\left(n-2\right)$, thus $C=2K+C''$. Therefore
\[
    \left.\frac{\partial G}{\partial\left\{x,y,w\right\}}\right|_K=
    \left.\frac{\partial \left.G\right|_E}
               {\partial\left\{x,y,w\right\}}\right|_K=0.
\]
Now $I=\left\{\partial G/\partial z|_E=0\right\}$ is a plane
curve of degree $n-1$, meeting all $2n$ nodes, all of whose
intersections with $K$ induce singularities on $Y$. So
$I.K\geq 2n>2\left(n-1\right)$, thus $I=K+I'$ which means that
\[
    \left.\frac{\partial G}{\partial z}\right|_K=
    \left.\left.\frac{\partial G}{\partial z}\right|_E\right|_K=0.
\]
So $Y$ is singular along $K$, this contradicts that
$K$ contains isolated singularities of $Y$.
\medskip\\
\indent So $p_0$ would induce a singular curve on $X_8$ which itself
would induce
at least one singularity on every plane $E_j$, $j=0,\ldots,n-1$ which is
not a  node, contradiction. Therefore the surface
\begin{eqnarray*}
    X_8 & = & \left\{64\left(x^2-w^2\right)\left(y^2-w^2\right)
              \left(\left(x+y\right)^2-2\,w^2\right)
              \left(\left(x-y\right)^2-2\,w^2\right)\right.\\
        &   & -\left[
                  -4\left(1+\sqrt{2}\right)\left(x^2+y^2\right)^2
                  +\left(
                      8\left(2+\sqrt{2}\right){z}^{2}
                      +2\left(2+7\,\sqrt{2}\right)w^2
                  \right )\left(x^2+y^2\right)
              \right.\\
        &   & \left.\left.
                  -16\,z^4
                  +8\left(1-2\,\sqrt{2}\right)z^2w^2
                  -\left(1+12\,\sqrt{2}\right)w^4\right]^2=0\right\}
\end{eqnarray*}
has exactly 168 nodes and no other singularities.


\begin{thebibliography}{[Ma]}
%
\bibitem[B]{barth}
    Barth,~W.:
    {\it Two Projective Surfaces with many Nodes, admitting the Symmetries
         of the Icosahedron.}
    (1995), to appear in Journal of Algebraic Geometry.
%
\bibitem[C]{cataneseceresa}
    Catanese,~F.; Ceresa,~G.:
    {\it Constructing Sextic Surfaces with a Given Number $d$ of Nodes.}
    Journal of Pure and Applied Algebra 23 (1982), 1--12.
%
\bibitem[G]{gallarati}
    Gallarati,~D.:
    {\it Una Superficie dell'ottavo Ordine con 160 Nodi.}
    Accad.\ Ligure Sci.\ Lett.\ 14 (1957), 1--7.
%
\bibitem[K]{kreiss}
    Kreiss,~H.~O.:
    {\it \"Uber syzygetische Fl\"achen.}
    Annali di Matematica (2) 41 (1955), 105--111.
%
\bibitem[Ma]{maple}
    {\it Maple V Release 3 Notes.}
    2nd edition (1994).
%
\bibitem[Mi]{miyaoka}
    Miyaoka,~Y.:
    {\it The Maximal Number of Quotient Singularities on Surfaces
         with Given Numerical Invariants.}
    Mathematische Annalen 268 (1984), 159--171.
%
%
%
\end{thebibliography}
\end{document}